\documentclass[a4paper,11pt]{amsart}
\usepackage{graphicx}
\usepackage{mathptmx}
\usepackage{amsmath}
\usepackage{amssymb}
\usepackage{enumitem}
\usepackage{xcolor}
\usepackage{xparse}
\NewDocumentCommand{\eulerian}{omm}
 {%
  \genfrac<>{0pt}{}{#2}{#3}%
  \IfValueT{#1}{_{\!#1}}%
 }

 \newmuskip\pFqmuskip

\newcommand*\pFq[6][8]{%
  \begingroup 
  \pFqmuskip=#1mu\relax
  \mathchardef\normalcomma=\mathcode`,
  \mathcode`\,=\string"8000
  \begingroup\lccode`\~=`\,
  \lowercase{\endgroup\let~}\pFqcomma
  {}_{#2}F_{#3}{\left(\genfrac..{0pt}{}{#4}{#5}\bigg|#6\right)}%
  \endgroup
}
\newcommand{\pFqcomma}{{\normalcomma}\mskip\pFqmuskip}

\newtheorem{theorem}{Theorem}

\newtheorem{proposition}[theorem]{Proposition}

\begin{document}

\title[Some identities involving degenerate Stirling numbers ]{Some identities involving degenerate Stirling numbers arising from normal ordering}

\author{Taekyun  Kim}
\address{Department of Mathematics, Kwangwoon University, Seoul 139-701, Republic of Korea}
\email{tkkim@kw.ac.kr}

\author{DAE SAN KIM}
\address{Department of Mathematics, Sogang University, Seoul 121-742, Republic of Korea}
\email{dskim@sogang.ac.kr}


\maketitle

\begin{abstract}
In this paper, we derive some identities and recurrence relations for the degenerate Stirling numbers of the first kind and of the second kind which are degenerate versions of the ordinary Stirling numbers of the first kind and of the second kind. They are deduced from the normal oderings of degenerate integral powers of the number operator and their inversions, certain relations of boson operators and from the recurrence relations of the Stirling numbers themselves. Here we note that, while the normal ordering of an integral power of the number operator is expressed with the help of the Stirling numbers of the second kind, that of a degenerate integral power of the number operator is represented by means of the degenerate Stirling numbers of the second kind.

\end{abstract}

\section{Introduction}

The Stirling number of the second $S_{2}(n,k)$ is the number of ways to partition a set of $n$ objects into $k$ nonempty subsets (see \eqref{6}). The (signed) Stirling number of the first kind $S_{1}(n,k)$ is defined in such a way that the number of permutations of $n$ elements having exactly $k$ cycles is the nonnegative integer $(-1)^{n-k}S_{1}(n,k)=|S_{1}(n,k)|$ (see \eqref{5}). The degenerate Stirling numbers of the second kind $S_{2,\lambda}(n,k)$ (see \eqref{9}) and of the first kind $S_{1,\lambda}(n,k)$ (see \eqref{7}) appear most naturally when we replace the powers of $x$ by the generalized falling factorial polynomials $(x)_{k,\lambda}$ in the defining equations (see \eqref{5}, \eqref{6}, \eqref{7}, \eqref{9}).\par

Carlitz initiated a study of degenerate versions of some special numbers and polynomials in [2], where the degenerate Bernoulli and Euler numbers were investigated. It is remarkable that in recent years intensive studies have been done for degenerate versions of quite a few special polynomials and numbers and have yielded many interesting results (see [6,7] and the references therein). They have been explored by various methods, including mathematical physics, combinatorial methods, generating functions, umbral calculus techniques, $p$-adic analysis, differential equations, special functions, probability theory and analytic number theory.
It turns out that the degenerate Stirling numbers play an important role in this exploration for degenerate versions of many special numbers and polynomials. The normal ordering of an integral power of the number operator $a^{+}a$ in terms of boson operators $a$ and $a^{+}$ can be written in the form
\begin{equation*}
(a^{+}a)^{k}=\sum_{l=0}^{k}S_{2}(k,l)(a^{+})^{l}a^{l}.
\end{equation*}
The normal ordering of the degenerate $k$th power of the number operator $a^{+}a$, namely $(a^{+}a)_{k,\lambda}$, in terms of boson operators $a, a^{+}$ can be written in the form
\begin{equation}
(a^{+}a)_{k,\lambda}=\sum_{l=0}^{k}S_{2,\lambda}(k,l)(a^{+})^{l}a^{l},\label{1}
\end{equation}
where the generalized falling factorials $(x)_{n,\lambda}$ are given by \eqref{1}.\par
By inversion, from \eqref{1} we obtain
\begin{equation}
(a^{+})^{k}a^{k}=\sum_{l=0}^{k}S_{1,\lambda}(k,l)(a^{+}a)_{l,\lambda}.\label{2}
\end{equation}
The aim of this paper is to derive some identities and recurrence relations for the degenerate Stirling numbers of the first kind and of the second kind by using the normal ordering in \eqref{1} and its inversion in \eqref{2}, certain relations of boson operators and the recurrence relations of the Stirling numbers themselves.
In more detail, our main results are as follows. Firstly, we derive some recurrence relation and identity involving the degenerate Stirling numbers of the second kind in Theorems 3 and 5 by using \eqref{1} and certain relations for boson operators. Secondly, we find some recurrence relation and identity involving the degenerate Stirling numbers of the first kind in Theorems 4 and 6 by exploiting \eqref{2} and certain relations for boson operators.
Thirdly, we investigate some recurrence relations for the degenerate Stirling numbers of the second kind in Theorems 7 and 9 by using the recurrence relation in
\eqref{15}, and those for the degenerate Stirling numbers of the first kind in Theorems 8 and 11 by using the recurrence relation in \eqref{8}. For the rest of this section, we recall the facts that are needed throughout this paper.

For any $\lambda\in\mathbb{R}$, the degenerate exponential functions are defined by
\begin{equation}
e_{\lambda}^{x}(t)=\sum_{k=0}^{\infty}(x)_{n,\lambda}\frac{t^{k}}{k!},\quad (\mathrm{see}\ [6,7]), \label{3}
\end{equation}
where
\begin{equation}
(x)_{0,\lambda}=1,\quad (x)_{k,\lambda}=x(x-\lambda)\cdots\big(x-(k-1)\lambda\big),\quad (k \ge 1).\label{4}
\end{equation}
When $x=1$, we use the notation $e_{\lambda}(t)=e_{\lambda}^{1}(t)$. It is well known that the Stirling numbers of the first kind are defined by
\begin{equation}
(x)_{n}=\sum_{k=0}^{n}S_{1}(n,k)x^{k},\quad (n\ge 0),\quad (\mathrm{see}\ [1,2,3,9,10]),\label{5}
\end{equation}
where $(x)_{0}=1,\ (x)_{n}=x(x-1)\cdots(x-n+1),\quad (n\ge 1)$. \par
As the inversion formula of \eqref{5}, the Stirling numbers of the second kind are defined by
\begin{equation}
x^{n}=\sum_{k=0}^{n}S_{2}(n,k)(x)_{k},\quad (n\ge 0).\label{6}
\end{equation}
In [6], the degenerate Stirling numbers of the first kind are defined by
\begin{equation}
(x)_{n}=\sum_{k=0}^{n}S_{1,\lambda}(n,k)(x)_{k,\lambda},\quad (n\ge 0).\label{7}	
\end{equation}
From \eqref{7}, we note that
\begin{equation}
S_{1,\lambda}(k+1,l)=S_{1,\lambda}(k,l-1)-(k-l\lambda)S_{1,\lambda}(k,l),\label{8}
\end{equation}
where $k,l\in\mathbb{N}$ with $k\ge l$. \par
As the inverse formula of \eqref{7}, the degenerate Stirling numbers of the second kind are defined by
\begin{equation}
(x)_{n,\lambda}=\sum_{k=0}^{n}S_{2,\lambda}(n,k)(x)_{k},\quad (n\ge 0),\quad (\mathrm{see}\ [6]).\label{9}
\end{equation}
Note that $\displaystyle\lim_{\lambda\rightarrow 0}S_{2,\lambda}(n,k)=S_{2}(n,k)\displaystyle$ and $\displaystyle\lim_{\lambda\rightarrow 0}S_{1,\lambda}(n,k)=S_{1}(n,k)\displaystyle $.\par
By \eqref{7} and \eqref{9}, we easily obtain the next proposition.
\begin{proposition}
The following orthogonality and inverse relations hold true.
\begin{align*}
&\sum_{k=l}^{n}S_{1,\lambda}(n,k)S_{2,\lambda}(k,l)=\delta_{n,l},
\,\, \sum_{k=l}^{n}S_{2,\lambda}(n,k)S_{1,\lambda}(k,l)=\delta_{n,l},\,\,(0 \le l \le n), \\
&a_{n}=\sum_{k=0}^{n}S_{2,\lambda}(n,k)b_{k} \Longleftrightarrow b_{n}=\sum_{k=0}^{n}S_{1,\lambda}(n,k)a_{k},\\
&a_{n}=\sum_{k=n}^{m}S_{2,\lambda}(k,n)b_{k} \Longleftrightarrow b_{n}=\sum_{k=n}^{m}S_{1,\lambda}(k,n)a_{k},
\end{align*}
where $\delta_{n,l}$ is the Kronecker's symbol. \\
\end{proposition}

Let $a$ and $a^{+}$ be the boson annihilation and creation operators satisfying the commutation relation
\begin{equation}
[a,a^{+}]=aa^{+}-a^{+}a=1,\quad (\mathrm{see}\ [4,5,8]).\label{10}	
\end{equation}
It has been known for some time that the normal ordering of $(a^{+}a)^{n}$ has a close relation to the Stirling numbers of the second kind ([4,5,8]).
Indeed, the normal ordering of an integral power of the number operator $a^{+}a$ in terms of boson operators $a$ and $a^{+}$, that satisfy the commutation $[a,a^{+}]=aa^{+}-a^{+}a=1$, can be written in the form
\begin{equation}
(a^{+}a)^{k}=\sum_{l=0}^{k}S_{2}(k,l)(a^{+})^{l}a^{l},\quad (\mathrm{see}\ [4,5,8]). \label{11}
\end{equation}
The number states $|m\rangle,\ m=1,2,\dots$, are defined as
\begin{equation}
a|m\rangle=\sqrt{m}|m-1\rangle,\quad a^{+}|m\rangle=\sqrt{m+1}|m+1\rangle.\label{12}
\end{equation}
By \eqref{12}, we get $a^{+}a|m\rangle=m|m\rangle$, (see [4,5,8]). The coherent states $|z\rangle$, where $z$ is a complex number, satisfy $a|z\rangle=z|z\rangle,$ $z|z\rangle=1$. To show a connection to coherent states, we recall that the harmonic oscillator has Hamiltonian $\hat{n}=a^{+}a$ (neglecting the zero point energy) and the usual eigenstates $|n\rangle$ (for $n\in\mathbb{N}$) satisfying $\hat{n}|n\rangle=n|n\rangle$ and $\langle m|n\rangle=\delta_{m,n}$, where $\delta_{m,n}$ is the Kronecker's symbol. In this paper, we derive some identities involving the degenerate Stirling numbers arising from the normal ordering of degenerate integral powers of the number operator $a^{+}a$ in terms of boson operators $a$ and $a^{+}$,
which are given by
\begin{displaymath}
(a^{+}a)_{k,\lambda}=\sum_{l=0}^{k}S_{2,\lambda}(k,l)(a^{+})^{l}a^{l}.
\end{displaymath}
In this paper, we study the degenerate Stirling numbers associated with the number operator $a^{+},a$ that satisfy $[a,a^{+}]=aa^{+}-a^{+}a=1$. Indeed, we give some new formulae of the degenerate Stirling numbers which are derived from an integral power of number operator $a^{+}a$ in terms of boson operators $a$ and $a^{+}$ with $$(a^{+}a)_{k,\lambda}=\sum_{l=0}^{k}S_{2,\lambda}(k,l)(a^{+})^{l}a^{l}.$$

\section{Some identities involving degenerate Stirling numbers}
We recall that the bosonic commutation relation
$[a,a^{+}]=aa^{+}-a^{+}a=1$ can be realized formally in a suitable space of functions by letting $a=\frac{d}{dx}$ and $a^{+}=x$ (the operator of multiplication by $x$). It is known that
\begin{equation}
\bigg(x\frac{d}{dx}\bigg)_{n,\lambda}f(x)=\sum_{k=0}^{n}S_{2,\lambda}(n,k)x^{k}\bigg(\frac{d}{dx}\bigg)^{k}f(x),\quad (\mathrm{see}\ [7]),\label{13}
\end{equation}
where $n$ is a nonnegative integer. \par
The equation \eqref{13} can be written as
\begin{equation}
(a^{+}a)_{n,\lambda}=\sum_{k=0}^{n}S_{2,\lambda}(n,k)(a^{+})^{k}a^{k}.\label{14}	
\end{equation}
From \eqref{9}, we note that
\begin{equation}
S_{2,\lambda}(k+1,l)=S_{2,\lambda}(k,l-1)+(l-k\lambda)S_{2,\lambda}(k,l),\label{15}
\end{equation}
where $k,l\in\mathbb{N}$ with $k\ge l$. It is easy to show that
\begin{displaymath}
	[a,\hat{n}]=a,\quad [\hat{n},a^{+}]=a^{+}.
\end{displaymath}
Inverting \eqref{13} by using Proposition 1, we have
\begin{equation}
x^{k}\bigg(\frac{d}{dx}\bigg)^{k}f(x)=\sum_{m=0}^{k}S_{1,\lambda}(k,m)\bigg(x\frac{d}{dx}\bigg)_{m,\lambda}f(x).\label{16}	
\end{equation}
In view of \eqref{16}, the normal ordering of a degenerate integral power of the number operator $a^{+}a$ in terms of boson operators $a,a^{+}$ can be rewritten in the form
\begin{align}
(a^{+})^{k}a^{k}&=\sum_{m=0}^{k}S_{1,\lambda}(k,m)(a^{+}a)_{m,\lambda}\label{17} \\
&=\sum_{m=0}^{k}S_{1,\lambda}(k,m)(\hat{n})_{m,\lambda}=(\hat{n})_{k},\nonumber	
\end{align}
where $k$ is a positive integer.
\begin{proposition}
For $k\in\mathbb{N}$, we have
\begin{align*}
(\hat{n})_{k}=(a^{+})^{k}a^{k}&=\sum_{m=0}^{k}S_{1,\lambda}(k,m)(a^{+}a)_{m,\lambda}
=\sum_{m=0}^{k}S_{1,\lambda}(k,m)(\hat{n})_{m,\lambda}.\end{align*}
\end{proposition}

We note that
\begin{align}
a^{+}(\hat{n}+1-\lambda)_{k,\lambda}a &=a^{+}\sum_{l=0}^{k}\binom{k}{l}(\hat{n})_{l,\lambda}(1-\lambda)_{k-l,\lambda}a \label{18}\\
&=\sum_{l=0}^{k}\binom{k}{l}(1-\lambda)_{k-l,\lambda}\sum_{m=0}^{l}S_{2,\lambda}(l,m)(a^{+})^{m+1}a^{m+1}\nonumber \\
&=\sum_{m=0}^{k}\bigg(\sum_{l=m}^{k}\binom{k}{l}(1-\lambda)_{k-l,\lambda}S_{2,\lambda}(l,m)\bigg)(a^{+})^{m+1}a^{m+1}. \nonumber	
\end{align}
By \eqref{10}, we also have
 \begin{align}
 (\hat{n})_{k+1,\lambda}&=(\hat{n}-\lambda)_{k,\lambda}\hat{n}=a^{+}\big((\hat{n}+1-\lambda)\cdots(\hat{n}+1-k\lambda)\big)a\label{19}\\
 &=a^{+}(\hat{n}+1-\lambda)_{k,\lambda}a.\nonumber
 \end{align}
By \eqref{19} and \eqref{14}, we get
\begin{align}
a^{+}(\hat{n}+1-\lambda)_{k,\lambda}a&=(\hat{n})_{k+1,\lambda}=(a^{+}a)_{k+1,\lambda}\nonumber \\
&=\sum_{m=0}^{k+1}S_{2,\lambda}(k+1,m)(a^{+})^{m}a^{m}\label{20} \\
&=\sum_{m=1}^{k+1}S_{2,\lambda}(k+1,m)(a^{+})^{m}a^{m}\nonumber \\
&=\sum_{m=0}^{k}S_{2,\lambda}(k+1,m+1)(a^{+})^{m+1}a^{m+1}.\nonumber
\end{align}
Therefore, by \eqref{18} and \eqref{20}, we obtain the following theorem.
\begin{theorem}
For $m,k\in\mathbb{Z}$ with $k\ge m \ge 0$, we have
\begin{displaymath}
\sum_{l=m}^{k}\binom{k}{l}(1-\lambda)_{k-l,\lambda}S_{2,\lambda}(l,m)=S_{2,\lambda}(k+1,m+1).
\end{displaymath}	
\end{theorem}
From Proposition 2, we note that
\begin{align}
(a^{+})^{k+1}a^{k+1}&=\sum_{m=0}^{k+1}S_{1,\lambda}(k+1,m)(a^{+}a)_{m,\lambda} \label{21}\\
&=\sum_{m=1}^{k+1}S_{1,\lambda}(k+1,m)(\hat{n})_{m,\lambda} \nonumber \\
&=\sum_{m=0}^{k}S_{1,\lambda}(k+1,m+1)(\hat{n})_{m+1,\lambda},\quad (k \ge 0).\nonumber
\end{align}
The degenerate rising factorial sequence is defined by
\begin{equation}
\langle x\rangle_{0,\lambda}=1,\quad \langle x\rangle_{n,\lambda}=x(x+\lambda)\cdots(x+(n-1)\lambda),\quad (n\ge 1).\label{22}	
\end{equation}
Now, from \eqref{10} we note that
\begin{equation}
a^{+}(\hat{n})_{k,\lambda}a=(\hat{n}-1)_{k,\lambda}\hat{n}=\hat{n}(\hat{n}-1)_{k,\lambda},\quad(k\in\mathbb{N}). \label{23}
\end{equation}
On the other hand, by \eqref{17} and \eqref{23}, we get
\begin{align}
&(a^{+})^{k+1}a^{k+1}=\sum_{l=0}^{k}S_{1,\lambda}(k,l)a^{+}(\hat{n})_{l,\lambda}a=\sum_{l=0}^{k}S_{1,\lambda}(k,l)(\hat{n}-1)_{l,\lambda}\hat{n} \label{24}\\
&=\sum_{l=0}^{k}S_{1,\lambda}(k,l)\sum_{m=0}^{l}\binom{l}{m}(-1)^{l-m}\langle 1\rangle_{l-m,\lambda}(\hat{n})_{m,\lambda}\big(\hat{n}-m\lambda+m\lambda) \nonumber \\
&=\sum_{l=0}^{k}S_{1,\lambda}(k,l)\sum_{m=0}^{l}\binom{l}{m}(-1)^{l-m}\langle 1\rangle_{l-m,\lambda}(\hat{n})_{m+1,\lambda}\nonumber \\
&\quad\quad +\lambda\sum_{l=0}^{k}S_{1,\lambda}(k,l)\sum_{m=1}^{l+1}\binom{l}{m}(-1)^{l-m}\langle 1\rangle_{l-m,\lambda}(\hat{n})_{m,\lambda}m\nonumber\\
&=\sum_{m=0}^{k}\sum_{l=m}^{k}S_{1,\lambda}(k,l)\binom{l}{m}(-1)^{l-m}\langle 1\rangle_{l-m,\lambda}(\hat{n})_{m+1,\lambda}\nonumber \\
&\quad\quad +\sum_{m=0}^{k}\sum_{l=m}^{k}S_{1,\lambda}(k,l)(m+1)\lambda\binom{l}{m+1}(-1)^{l-m-1}\langle 1\rangle_{l-m-1,\lambda}(\hat{n})_{m+1,\lambda}.\nonumber
\end{align}
Therefore, by \eqref{21} and \eqref{24}, we obtain the following theorem.
\begin{theorem}
For $m,k\in\mathbb{Z}$ with $k \ge m \ge 0$, we have
\begin{align*}
&S_{1,\lambda}(k+1,m+1) \\
&\quad =\sum_{l=m}^{k}S_{1,\lambda}(k,l)\bigg\{\binom{l}{m}\langle 1\rangle_{l-m,\lambda}(-1)
^{l-m}+l\lambda\binom{l-1}{m}(-1)^{l-m-1}\langle 1\rangle_{l-m-1,\lambda}\bigg\}.
\end{align*}
\end{theorem}
From \eqref{23}, we get
\begin{align}
&a^{+}(\hat{n})_{k,\lambda}a=\hat{n}(\hat{n}-1)_{k,\lambda}=\sum_{m=0}^{k}\binom{k}{m}(-1)^{k-m}\langle 1\rangle_{k-m,\lambda}\hat{n}(\hat{n})_{m,\lambda}\label{25} \\
&=\sum_{m=0}^{k}\binom{k}{m}(-1)^{k-m}\langle 1\rangle_{k-m,\lambda}(\hat{n}-m\lambda+m\lambda)(\hat{n})_{m,\lambda}\nonumber \\
&=\sum_{m=0}^{k}\binom{k}{m}(-1)^{k-m}\langle 1\rangle_{k-m,\lambda}\big((\hat{n})_{m+1,\lambda}+m\lambda(\hat{n})_{m,\lambda}\big) \nonumber \\
&= \sum_{m=0}^{k}\binom{k}{m}(-1)^{k-m}\langle 1\rangle_{k-m,\lambda}(\hat{n})_{m+1,\lambda}\nonumber\\
& \quad +\sum_{m=0}^{k}\binom{k}{m+1}(-1)^{k-m-1}\langle 1\rangle_{k-m-1,\lambda}(m+1)\lambda (\hat{n})_{m+1,\lambda} \nonumber\\
&=\sum_{m=0}^{k}\binom{k}{m}(-1)^{k-m}\langle 1\rangle_{k-m,\lambda}\sum_{p=1}^{m+1}S_{2,\lambda}(m+1,p)(a^{+})^{p}a^{p} \nonumber \\
&\quad +\lambda\sum_{m=0}^{k}\binom{k}{m+1}(-1)^{k-m-1}\langle 1\rangle_{k-m-1,\lambda}(m+1)\sum_{p=1}^{m+1}S_{2,\lambda}(m+1,p)(a^{+})^{p}a^{p}. \nonumber \\
&=\sum_{p=1}^{k+1}\bigg(\sum_{m=p-1}^{k}\binom{k}{m}(-1)^{k-m}\langle 1\rangle_{k-m,\lambda}S_{2,\lambda}(m+1,p)\bigg)(a^{+})^{p}a^{p} \nonumber \\
&\quad +\lambda\sum_{p=1}^{k+1}\sum_{m=p-1}^{k}\binom{k}{m+1}(-1)^{k-m-1}\langle 1\rangle_{k-m-1,\lambda}(m+1)S_{2,\lambda}(m+1,p)(a^{+})^{p}a^{p}\nonumber \\
&=\sum_{p=0}^{k}\bigg\{\sum_{m=p}^{k}\binom{k}{m}(-1)^{k-m}\langle 1\rangle_{k-m,\lambda}S_{2,\lambda}(m+1,p+1)\bigg\}(a^{+})^{p+1}a^{p+1} \nonumber \\
&\quad +\lambda\sum_{p=0}^{k}\sum_{m=p}^{k}k\binom{k-1}{m}(-1)^{k-m-1}\langle 1\rangle_{k-m-1,\lambda}S_{2,\lambda}(m+1,p+1)(a^{+})^{p+1}a^{p+1}. \nonumber
\end{align}
We observe that
\begin{align}
a^{+}(\hat{n})_{k,\lambda}a &= a^{+}\bigg(\sum_{p=0}^{k}S_{2,\lambda}(k,p)(a^{+})^{p}a^{p}\bigg)a \label{26}\\
&=\sum_{p=0}^{k}S_{2,\lambda}(k,p)(a^{+})^{p+1}a^{p+1}.\nonumber	
\end{align}
Therefore, by \eqref{25} and \eqref{26}, we obtain the following theorem.
\begin{theorem}
For $p,k\in\mathbb{Z}$ with $0\le p \le k$, we have
\begin{align*}
	S_{2,\lambda}(k,p)=&\sum_{m=p}^{k}\bigg\{\binom{k}{m}(-1)^{k-m}\langle 1\rangle_{k-m,\lambda}S_{2,\lambda}(m+1,p+1)\\
	&\quad +\lambda k\binom{k-1}{m}(-1)^{k-m-1}\langle 1\rangle_{k-m-1,\lambda}S_{2,\lambda}(m+1,p+1)\bigg\}.
\end{align*}	
\end{theorem}

By \eqref{17}, we obtain
\begin{align}
&(a^{+})^{k+1}a^{k+1}=\sum_{p=0}^{k+1}S_{1,\lambda}(k+1,p)(a^{+}a)_{p,\lambda}=\sum_{p=1}^{k+1}S_{1,\lambda}(k+1,p)(\hat{n})_{p,\lambda} \label{27} \\
&=\sum_{p=1}^{k+1}S_{1,\lambda}(k+1,p)\hat{n}(\hat{n}-1+1-\lambda)_{p-1,\lambda}\nonumber\\
&=\sum_{p=1}^{k+1}S_{1,\lambda}(k+1,p)\hat{n}\sum_{l=0}^{p-1}\binom{p-1}{l}(1-\lambda)_{p-1-l,\lambda}(\hat{n}-1)_{l,\lambda}\nonumber\\
&=\sum_{l=0}^{k}\bigg(\sum_{p=l+1}^{k+1}S_{1,\lambda}(k+1,p)\binom{p-1}{l}(1-\lambda)_{p-1-l,\lambda}\bigg)\hat{n}(\hat{n}-1)_{l,\lambda}\nonumber\\
&=\sum_{l=0}^{k}\bigg(\sum_{p=l}^{k}S_{1,\lambda}(k+1,p+1)\binom{p}{l}(1-\lambda)_{p-l,\lambda}\bigg)\hat{n}(\hat{n}-1)_{l,\lambda}.\nonumber
\end{align}
Therefore, from the first line of \eqref{24} and \eqref{27}, we obtain the following theorem.
\begin{theorem}
For $k,l\in\mathbb{Z}$ with $k\ge l \ge 0$, we have
\begin{displaymath}
S_{1,\lambda}(k,l)=\sum_{p=l}^{k}S_{1,\lambda}(k+1,p+1)\binom{p}{l}(1-\lambda)_{p-l,\lambda}.
\end{displaymath}	
\end{theorem}
From \eqref{9}, we note that
\begin{equation}
\frac{1}{k}\big(e_{\lambda}(t)-1\big)^{k}=\sum_{n=k}^{\infty}S_{2,\lambda}(n,k)\frac{t^{n}}{n!},\quad (k\ge 0).\label{28}
\end{equation}
Thus, by \eqref{28}, we get
\begin{align}
\sum_{n=k}^{\infty}S_{2,\lambda}(n,k)\frac{t^{n}}{n!}&=\frac{1}{k!}\big(e_{\lambda}(t)-1\big)^{k}=\frac{1}{k!}\sum_{m=0}^{k}\binom{k}{m}(-1)^{k-m}e_{\lambda}^{m}(t)\label{29} \\
&=\sum_{n=0}^{\infty}\bigg(\frac{1}{k!}\sum_{m=0}^{k}\binom{k}{m}(-1)^{k-m}(m)_{n,\lambda}\bigg)\frac{t^{n}}{n!}.\nonumber
\end{align}
Comparing the coefficients on both sides of \eqref{29}, we have
\begin{displaymath}
\sum_{m=0}^{k}\binom{k}{m}(-1)^{k-m}(m)_{n,\lambda}=\left\{\begin{array}{ccc}
k!S_{2,\lambda}(n,k), & \textrm{if $n\ge k,$}\\
0, & \textrm{if $0 \le n<k$}.
\end{array}\right.
\end{displaymath}
By \eqref{15}, we get
\begin{align}
&S_{2,\lambda}(k+1,l+1)=S_{2,\lambda}(k,l)+(l+1-k\lambda)S_{2,\lambda}(k,l+1) \label{30}\\
&=S_{2,\lambda}(k,l)+(l+1-k\lambda)\nonumber \\
&\quad\quad \times \big\{S_{2,\lambda}(k-1,l)+(l+1-(k-1)\lambda)S_{2,\lambda}(k-1,l+1)\big\}\nonumber \\
&=S_{2,\lambda}(k,l)+(l+1-k\lambda)S_{2,\lambda}(k-1,l) \nonumber \\
&\quad\quad +\langle l+1-k\lambda \rangle_{2,\lambda}S_{2,\lambda}(k-1,l+1) \nonumber \\
&=S_{2,\lambda}(k,l)+(l+1-k\lambda)S_{2,\lambda}(k-1,l)+\langle l+1-k\lambda \rangle_{2,\lambda}S_{2,\lambda}(k-2,l)\nonumber\\
&\quad+\cdots + \langle l+1-k\lambda \rangle_{k-l,\lambda}S_{2,\lambda}(k,l) \nonumber \\
&=\sum_{m=l}^{k} \langle l+1-k\lambda \rangle_{k-m,\lambda}S_{2,\lambda}(l,l).\nonumber
\end{align}
Therefore, by \eqref{30}, we obtain the following theorem.
\begin{theorem}
For $l,k\in\mathbb{Z}$ with $0\le l\le k$, we have
\begin{displaymath}
S_{2,\lambda}(k+1,l+1)=\sum_{m=l}^{k}\langle l+1-k\lambda \rangle_{k-m,\lambda}S_{2,\lambda}(m,l).
\end{displaymath}
\end{theorem}
From \eqref{8}, we note that
\begin{align}
&S_{1,\lambda}(k+1,m+1)=S_{1,\lambda}(k,m)-\big(k-(m+1)\lambda\big)S_{1,\lambda}(k,m+1)\label{31}\\
&=S_{1,\lambda}(k,m)-\big(k-(m+1)\lambda\big)\nonumber\\
&\quad\quad \times\Big(S_{1,\lambda}(k-1,m)-((k-1)-(m+1)\lambda)S_{1,\lambda}(k-1,m+1)\Big)\nonumber \\
&=S_{1,\lambda}(k,m)-(k-(m+1)\lambda)S_{1,\lambda}(k-1,m)\nonumber \\
&\quad\quad +(-1)^{2}(k-(m+1)\lambda)_{2} S_{1,\lambda}(k-1,m+1)\nonumber \\
&=S_{1,\lambda}(k,m)-(k-(m+1)\lambda)S_{1,\lambda}(k-1,m) \nonumber \\
&\quad\quad +(-1)^{2}(k-(m+1)\lambda)_{2} S_{1,\lambda}(k-1,m+1)\nonumber \\
&\quad\quad +\cdots+(-1)^{k-m}(k-(m+1)\lambda)_{k-m}S_{1,\lambda}(m,m)\nonumber \\
&=\sum_{l=m}^{k}(-1)^{k-l}\big(k-(m+1)\lambda)_{k-l}S_{1,\lambda}(l,m).\nonumber
\end{align}
Therefore, by \eqref{31}, we obtain the following theorem.
\begin{theorem}
For $m,k\in\mathbb{Z}$ with $k \ge  m \ge 0$, we have
\begin{displaymath}
S_{1,\lambda}(k+1,m+1)=\sum_{l=m}^{k}(-1)^{k-l}\big(k-(m+1)\lambda\big)_{k-l}S_{1,\lambda}(l,m).
\end{displaymath}	
\end{theorem}
By \eqref{15}, we get
\begin{align}
&S_{2,\lambda}(m+k+1,m)=S_{2,\lambda}(m+k,m-1)+\big(m-(m+k)\lambda\big)S_{2,\lambda}(m+k,m)\label{32}\\
&=S_{2,\lambda}(m+k-1,m-2)+\big((m-1)-(m+k-1)\lambda\big)S_{2,\lambda}(m+k-1,m-1) \nonumber \\
&\quad\quad +\big(m-(m+k)\lambda\big)S_{2,\lambda}(m+k,m)\nonumber \\
&=\big(m-(m+k)\lambda\big)S_{2,\lambda}(m+k,m)+\big((m-1)-(m+k-1)\lambda\big)\nonumber \\
&\quad\quad \times S_{2,\lambda}(m+k-1,m-1)+\cdots+(0-k\lambda)S_{2,\lambda}(k,0)	\nonumber \\
&=\sum_{l=0}^{m}\big(l-(k+l)\lambda\big)S_{2,\lambda}(k+l,l).\nonumber
\end{align}
Therefore, by \eqref{32}, we obtain the following theorem.
\begin{theorem}
For $k,m\in\mathbb{Z}$ with $k,m\ge 0$, we have
\begin{displaymath}
	S_{2,\lambda}(m+k+1,m)=\sum_{l=0}^{m}\big(l-(k+l)\lambda\big)S_{2,\lambda}(k+l,l).
\end{displaymath}	
\end{theorem}
From Theorem 3, we note that
\begin{displaymath}
	S_{2,\lambda}(k+1,m+1)=\sum_{l=m}^{k}\binom{k}{l}(1-\lambda)_{k-l,\lambda}S_{2,\lambda}(l,m).
\end{displaymath}
Now, by using Propositon 1 we have
\begin{align}
\sum_{l=m}^{k}S_{2,\lambda}(k+1,l+1)S_{1,\lambda}(l,m)&=\sum_{l=m}^{k}\sum_{p=l}^{k}\binom{k}{p}(1-\lambda)_{k-p,\lambda}S_{2,\lambda}(p,l)S_{1,\lambda}(l,m).\label{33}\\
&=\sum_{p=m}^{k}\binom{k}{p}(1-\lambda)_{k-p,\lambda}\sum_{l=m}^{p}S_{2,\lambda}(p,l)S_{1,\lambda}(l,m)\nonumber \\
&=\binom{k}{m}(1-\lambda)_{k-m,\lambda}.\nonumber	
\end{align}
Therefore, by \eqref{33}, we obtain the following theorem.
\begin{theorem}
For $m,k\in\mathbb{Z}$ with $k \ge m \ge 0$, we have
\begin{displaymath}
	\binom{k}{m}=\frac{1}{(1-\lambda)_{k-m,\lambda}}\sum_{l=m}^{k}S_{2,\lambda}(k+1,l+1)S_{1,\lambda}(l,m).
\end{displaymath}	
\end{theorem}
By \eqref{8}, we get
\begin{align}
&S_{1,\lambda}(m+k+1,m)=S_{1,\lambda}(m+k,m-1)-(m+k-m\lambda)S_{1,\lambda}(m+k,m)\label{34} \\
&=S_{1,\lambda}(m+k-1,m-2)-(m+k-1-(m-1)\lambda))S_{1,\lambda}(m+k-1,m-1)\nonumber\\
 &\quad\quad -(m+k-m\lambda)S_{1,\lambda}(m+k,m)\nonumber \\
 &=-(m+k-m\lambda)S_{1,\lambda}(m+k,m)-(m+k-1-(m-1)\lambda)\nonumber \\
 &\quad\quad \times S_{1,\lambda}(m+k-1,m-1)- \cdots-k S_{1,\lambda}(k,0)\nonumber\\
 &=-\sum_{l=0}^{m}(k+l-l\lambda)S_{1,\lambda}(k+l,l).\nonumber
 \end{align}
Therefore, by \eqref{34}, we obtain the following theorem.
\begin{theorem}
For $m,k\in\mathbb{Z}$ with $m,k\ge 0$, we have
\begin{displaymath}
S_{1,\lambda}(m+k+1,m)=-\sum_{l=0}^{m}(k+l-l\lambda)S_{1,\lambda}(k+l,l).
\end{displaymath}
\end{theorem}
\section{Conclusion}

In recent years, studying degenerate versions of some special numbers and polynomials have drawn the attention of many mathematicians with their regained interests not only in combinatorial and arithmetical properties but also in applications to differential equations, identities of symmetry and probability theory. These degenerate versions include the degenerate Stirling numbers of the first and second kinds, degenerate Bernoulli numbers of the second kind and degenerate Bell numbers and polynomials. \par
As a degenerate version of the well known normal ordering of an integral power of the number operator, we considered the normal ordering of a degenerate integral power of the number operator in terms of boson operators and its inversion as well.
We derived some identities and recurrence relations for the degenerate Stirling numbers of the first kind and of the second kind by using the normal ordering in \eqref{1} and its inversion in \eqref{2}, certain relations of boson operators and the recurrence relations of the Stirling numbers themselves.  \par
It is one of our future projects to continue to explore various degenerate versions of many special polynomials and numbers by using various methods mentioned in the Introduction.

\vspace{ 0.5cm}

{ \bf Conflict of Interest}

We have no conflicts of interest to disclose.

\vspace{ 0.5cm}
{\bf DATA AVAILABILITY}

Data sharing is not applicable to this article as no new data were created or analyzed in this study.

\end{document}